# Theodorus' lesson in Plato's *Theaetetus* (147d1-d6) Revisited.

# A New Perspective

**Luc Brisson**
Centre Jean Pépin
CNRS-UMR 8230

**Salomon Ofman**
CNRS-Institut mathématique de Jussieu-
Paris Rive Gauche
Histoire des Sciences mathématiques

**Abstract.** This article is the first part of a study of the so-called 'mathematical part' of Plato's *Theaetetus* (147d-148b). The subject of this 'mathematical part' is the irrationality, one of the most important topics in early Greek mathematics. As of huge interest for mathematicians, historians of mathematics as well as of philosophy, there had been an avalanche of studies about it. In our work, we revisit this question, for we think something is missing: a global analysis of Plato's text, from these three points of view simultaneously: history, mathematics and philosophy. It is what we have undertook through a new translation, a new interpretation of the mathematical lesson about irrational magnitudes and a novel interpretation of the whole passage from these three points of view.

Our guideline is considering Plato's writings seriously, not as some playful work. This simple rule is indeed surprisingly constraining, but it allows us to get a rare direct glance inside pre-Euclidean mathematics, in contradiction with the 'Main Standard Interpretation' prevailing in history of mathematics as well as in history of philosophy.

This study had been divided in two parts for editorial reasons. In the present article, we propose an analysis of the first part of this 'mathematical part', Theodorus' lesson. In the second article (Brisson-Ofman (to appear)), we present the sequel of the lesson and a philosophical interpretation of the 'mathematical part' within the framework of the entire dialogue.

Both articles form a whole. They are both aimed to an audience without any particular mathematical background, and require only elementary mathematical knowledge, essentially of high school-level. Some more delicate points are nevertheless developed in the **Appendices**.

**Résumé.** *Cet article est la première partie d'une étude de ce qu'on appelle la 'partie mathématique' du* Théétète *de Platon (147d-148b). Le sujet de cette 'partie mathématique' est l'irrationalité qui, en tant qu'il concerne aussi bien les mathématiciens, les historiens des mathématiques que ceux de la philosophie, a été l'objet d'une avalanche d'études et de commentaires. Nous reprenons néanmoins cette question, car nous pensons qu'une approche globale fait défaut : une analyse du texte, simultanément de ces trois points de vue : histoire, mathématiques et philosophie. C'est ce que nous avons entrepris ici, en proposant une* *nouvelle traduction, une nouvelle interprétation de la leçon mathématique de Théodore sur les grandeurs irrationnelles et une nouvelle interprétation du passage tout entier.*

*Notre fil conducteur constant a été de prendre le texte de Platon au sérieux, non pas comme une sorte de récit fantaisiste. Cette simple règle est certes étonnamment contraignante, mais elle permet d'obtenir un aperçu direct sur les mathématiques pré-euclidiennes. Il est en contradiction ouverte avec 'l'interprétation standard principale' prévalant en histoire des mathématiques aussi bien qu'en histoire de la philosophie.*

*Pour des raisons éditoriales, nous avons divisé ce travail en deux articles. Dans celui-ci, nous étudions la première partie du passage, la leçon de Théodore. Dans le second (Brisson-Ofman (to appear)), nous présentons la suite et la fin du passage, ainsi qu'une interprétation philosophique globale de la 'partie mathématique, dans le cadre du dialogue dans sa totalité.*

*Ces deux articles forment un tout. Ils se destinent tous deux à un public sans formation mathématique particulière, leur compréhension ne supposant que des connaissances mathématiques très élémentaires, essentiellement celles des premières années de collège. Pour le lecteur curieux, certains points plus délicats sont néanmoins développés dans les* **Annexes**.

## I. Introduction

The so-called 'mathematical part' of the *Theaetetus* (147d-148b) is divided in two parts. The first, presented in this article, concerns the beginning of the passage (147d1-d6). The second part is studied in its sequel (Brisson-Ofman (to appear)). Our fundamental guideline is to stick closely to Plato's text, in other words, to take it seriously and thus, to avoid any use of grand hypotheses about Plato's supposed purposes.

As a matter of fact, it is unreasonable to assert that some passages are fictitious, when Plato presents them as realistic. His works are indeed filled with extensive and famous myths, as well as with imaginary scenes like dreams. In each case, however, Plato makes clear the reasons for introducing these fictions.

The present article challenges the view of the 'Main Standard Interpretations' (or MSI) [1] according to which the modern commentator is free to decide what in Plato's text is true or not, for the following reasons:

- Plato had a huge amount of knowledge on almost any scientific field of his time, whether directly or through scientists attending the Academy. This is obvious to anyone who has read the *Statesman,* the *Timaeus* or the *Laws*.
- The Academy was open to people who excelled in many different fields, especially mathematics. And for instance, mathematicians would certainly have protested if Plato had written absurd things about their teaching, learning or work.
- Moreover, it would certainly not have convinced his reader or auditor, who in a small city like Athens knew exactly the facts recounted by Plato and were able to criticize him. Thus, false or dubious information concerning Socrates, Theodorus

- and Theaetetus would have made Plato appear ridiculous in their view, and his adversaries would have used it against him.
- Lastly, we must remember that the relations between the different schools of philosophy were not as polite as they are today, as we can see from the anecdote about Diogenes the Cynic throwing a plucked chicken during one of Plato's classes, because he did not like a definition given by the latter. [2] Against such a background, it is not certain Plato (or anyone else) would run the risk of claiming a total absurdity or of presenting some fictional character in one of the main role and in mixing him with real personages. [3]

We do not mean Plato describes real events as they happen. He is not a journalist. But he stages all his dialogues only from inescapable facts. [4] When one takes Plato's text seriously, one is faced by two difficult tasks that the standard interpreter may disregard. [5] The first is to show that such a literal interpretation is possible; the second that it is more comprehensive than the standard one. This is also one of the reasons for the length of this study, which we had to cut into two parts.

When a historian of philosophy talks about the Theaetetus today, he is immediately referred to the "modern standard interpretation", which is situated in the philosophical context of "logical empiricism". It differs from the interpretation that has been traditional from the Platonists of the anonymous commentator on the *Theaetetus* [6] to F. M. Cornford [7], H. Cherniss [8], L. Robin [9] and A. Diès [10]. This interpretation considers that if the search for a definition of science, which is the subject of the dialogue, ends up in a defeat, this is because its goal is 'to highlight the ruinous consequences of attempting epistemology without the backing of a Platonic metaphysics.' [11]

On the contrary, the 'largely dominant [interpretation], especially in English-language studies of Plato, for more than three decades', proposes an approach to the *Theaetetus* inspired in particular by the works of such Oxford scholars as Richard Robinson and G.E.L. Owen in the 1950s'.[12] One of the most important authors in this interpretive tendency is M. Burnyeat. [13] In an article published in 1978, [14] he tries to account for maieutics by leaving aside such metaphysical implications as the Forms and the soul separated from the body, and reducing it to an educational method intended to promote self-knowledge. More than 20 years later, he wrote a long preface to the translation by M.J. Levett, in which, without repudiating his 1978 article, he developed a new interpretation of the dialogue: "On Burnyeat's view, the *Theaetetus* is a dialectical exercise rather than a doctrinal one." [15] In order to begin a dialogue with the reader, Burnyeat, in his introduction/commentary, draws an opposition between a reading A, corresponding to the first type of interpretation, and a reading B, which corresponds to the second type, and indicates that he prefers the second one, inaugurated by Ryle and Owen.

The reading in question is an "analytic" one, which takes into consideration only the central part of the dialogue, in which we find the three definitions of science (*epistêmê*) [16] proposed from 151d-209a, but leaving aside the mathematical passage (147b-148d), as well as the passage on maieutics (*Theaetetus* 148e-151d) and the digression (172c3-177c5), which provide essential information for producing a philosophical interpretation that takes the intelligible forms and the soul into consideration. The analytical approach is based on a series of maneuvers intended to isolate a more or less lengthy passage of a dialogue. [17] This is the type of approach we will oppose here, basing ourselves on the mathematical passage. [18]

In parallel, for modern historians of mathematics, an author such as Plato – who, to top it all off, a philosopher – is scarcely credible. Thus, to establish the real existence of Theaetetus, Wilbur Knorr finds it necessary to call upon, rather than Plato, a text attributed (no doubt correctly) to Pappus, [19] of which only a late Arabic translation exists. [20] However, Pappus was a mathematician who lived a millennium after these events, and his interest in the history of his predecessors is not established. In this text, Eudemus of Rhodes is cited in support of the statements about Theaetetus, who concerns him only to emphasize Euclid's results. [21] What is more, quotations by authors of Late Antiquity are known for being subject to caution. They are often made from summaries compiled at second or third hand, since consulting papyrus rolls was anything but easy. In addition, with some exceptions, the quotation properly so called is mixed inextricably with the author's opinion, so this it is generally very delicate to separate one from the other. Finally, in the present case, Eudemus, the student of Aristotle, did indeed write a history of mathematics, but he could not be said to be a historian in the contemporary sense of the word, in search of sources and archives. Absent precise information on Theaetetus, who certainly died around 390 BCE, and hence very young (cf. *infra*, note 22), he probably was content to use the Platonic text as his source. This is the only text that provides direct information on Theaetetus, since the few other texts from Antiquity that discuss him are much later and may all be based on the eponymous dialogue.

We do not, of course, place in doubt the reality of a mathematician named Theaetetus. For the reasons cited above, that would be simply absurd. Nor is there any reason to suspect the information given by Pappus, taken from Eudemus even indirectly. The very existence of a long tradition, in a mathematical environment that showed little interest in individual history, makes them quite plausible. Nevertheless, we wish to emphasize here the incoherence of opposing a hypercritical skepticism against certain testimonies, in order to affirm, without solid foundations, that certain events are perfectly well established, especially when they comfort some thesis. We wish to react against this type of interpretation, restoring a central place in the history of philosophy and mathematics to the testimony of Plato. This is why it is necessary to submit the "axiomatic" hypotheses of the MSI to close criticism. [22]

In the article to follow this one, we would like to show that each of the three moments of the mathematical passage corresponds to the three definitions of knowledge proposed in the following sections of the text. The lesson of Theodorus, who constructs figures (147d1-6), corresponds to the definition of knowledge by sense perception (151d-171e, 177d-186e). The association of a type of figure with a type of number (147d6-148a4) corresponds to the definition of knowledge by true judgment (187a-200d), and the attempt to define power (148a4-148b3) corresponds to the definition of knowledge by true opinion accompanied by a *logos* (200d-210a).

This synthetic interpretation implies that one takes seriously the passage on maieutics, and the digression, and hence that one maintains the hypothesis of Forms grasped by the soul in the context of a process of reminiscence. [23] Socrates does not transmit any teaching to Theaetetus; he leads him to rid himself of his false opinions, with a view to disposing him to rediscover the true knowledge that his soul has contemplated before it became inscribed within a body. It then becomes imperative to recognize also the importance Socrates gives the next day to the continuation of the discussion in the *Sophist*.

The fundamental role of the mathematical part for the interpretation of the whole dialogue has been recognized since Antiquity. We will avoid the use of modern mathematical notions, except when explaining some important point as in the **Appendices**. Modern mathematics tends toward a universal algebraic language, while Greek mathematics was essentially geometrical; [24] hence the importance of drawings, which were sometimes even confused with the proof itself. [25] For instance, the modern product of two lengths A×B is represented as a rectangle of sides A and B (see for instance, *Theaetetus*, 148a). The reader must always keep in mind that the proof is carried out by means of a diagram, not through some algorithm, as algebraic notations would lead one to believe. Theodorus emphasizes this himself when he claims he was fed up with 'discourses' ('λόγοι') and had been leaving them to geometry. [26]

However, contrary to a recent trend in history of mathematics, we will certainly not consider ourselves bound by the language of Euclid's *Elements*. First, it was written more than a century after the dialogue, a period when mathematics underwent many deep changes. Looking in Euclid's *Elements* to find some answer to problems at the end of the $5^{th}$ century is like searching in Bourbaki to find demonstrations about mathematics questions in the beginning of the $19^{th}$ century. The best hints about the kind and the style of mathematics in Socrates' time are indeed found in Plato's texts, and to some extent, in some books of Aristotle. [27] Moreover, contrary to Platonic texts, the Euclidean books were strongly edited in the next centuries, in particular, but not only, by Theon of Alexandria. Thus we will not try to mimic the wording in the *Elements*, as often done in papers on early Greek mathematics, looking for some supposed proofs of pre-Euclidean results. [28] Almost nothing is known about the *forms* of pre-Euclidean mathematics. [29] We will always consider only demonstrations consistent with Plato's text, certainly known in the middle of the $5^{th}$ century and even much earlier (for instance by Babylonian and Egyptian mathematicians) that are moreover easily translated into Euclidean language. It is then useless to complicate matters for a modern reader by sticking to the often elliptical style of the *Elements*. For instance, we will freely speak of the 'area' of a square or a rectangle, when Euclid, and, in this particular case, probably the ancient Greek geometers, would have simply said 'square' or 'rectangle'.

### II. Staging

The dialogue opens with a magnificent prologue. We are at Megara, a city situated at the eastern end of the isthmus of Corinth, halfway between Corinth and Athens. The first character in the prologue is Euclid [30] of Megara. He is a member of Socrates' circle, and was later considered as the founder of the Megarian school. He meets Terpsion, who has been looking for him. Euclid explains that he was not in Megara, for he has just returned from Erinos, a city not far from Eleusis (20 kilometers west of Athens). This is where he has just accompanied Theaetetus, who has been wounded at a battle not far from Corinth [31] and is suffering from dysentery. Having arrived by boat at Megara, Theaetetus, on the verge of death, had insisted on coming home. On the way back from Erinos to Megara, Euclid recalls a discussion between Theaetetus and Socrates, just before the latter's condemnation and death. [32] Euclid, who was not present at this discussion, had written down what Socrates had recounted to him, and one of his slaves reads it out to him and to Terpsion.

Terpsion transports us from Megara to Athens in 399, where Socrates meets Theodorus, a geometer from Cyrene, [33] and questions him about the students who are following his teaching. Theodorus mentions the most brilliant of them, claiming that his physical appearance is as unprepossessing as Socrates', but that his psychic qualities are exceptional; and he introduces him to Socrates as Theaetetus.

Thus begins a dialogue between the three characters. Socrates soon raises doubts about Theodorus' expertise on physical beauty and ugliness, and Theaetetus is compelled, despite his reluctance, to agree with him. Then, as he does with Charmides, [34] Socrates changes registers. From Theaetetus' body, he switches to his soul, which Theodorus has just praised. Socrates wants to verify the validity of this praise, by enquiring into the knowledge that Theaetetus is trying to acquire, and more generally into science. Here begins a discussion about a definition (*logos*) of what science is. Asked to be the respondent, Theodorus refuses, and proposes Theaetetus instead.

The latter accepts tentatively, and Socrates says later, he will try, through the discussion, to help to make Theaetetus' soul give birth to a definition. [35] In response to Socrates' question about science, he gives a list of sciences: geometry, astronomy, (musical) harmony, calculation ('*logismos*'), all sciences already mentioned by Socrates, with the addition of shoemaking and other practical techniques ('*technai*')[36] (146c-d). This amalgam between techniques and sciences is already an example of the main subject, incommensurability. In the introduction of Plato's *Statesman*, where the same characters appear as in another dialogue of Plato, the *Sophist*, supposed to take place the day after the discussion recounted in *Theaetetus*, Socrates severely berates Theodorus when he puts politics and philosophy on the same level. And indeed, the latter has to admit there is no proportion between them (257a9-b8). This also shows the ambiguity of the term '*epistêmê*' which throughout the dialogue means mostly 'science' or 'scientific knowledge' for Socrates while for Theaetetus it means 'knowledge'.

Socrates reproaches him for his 'generosity' (146d4-5), previously praised by Theodorus as one of his qualities (144d3). Indeed, the young boy gave many examples, when only one definition ('*logos*') of science/knowledge had been sought. Then, suddenly, Theaetetus says that Socrates' explanations of his question 'what is science?', remind him of a similar problem. As he and a friend of his, Socrates' namesake, were thinking about Theodorus' recent lesson on powers, they raised a problem for themselves in connection with this lesson.

To distinguish between both 'Socrates', we will write Theaetetus' friend in italics (as '*Socrates*').

### III. The text

In order to avoid a circular interpretation, contrary to most translations, ours will closely follow the Greek, without striving for elegance or trying to paraphrase it. However, for a modern reader, it may seem cryptic. Thus, we will explain the main points in notes following it.

| THEAETETUS | ΘΕΑΙ. |
|---|---|
| **[147d]** Theodorus was drawing for us something about powers to prove that the ones of three 'feet' and of five 'feet' are not commensurable as length with one-foot long. Then he went on in this way, taking each case in turn, till the one of seventeen 'feet'. There for some reason he stopped. | Περὶ δυνάμεών τι ἡμῖν Θεόδωρος ὅδε ἔγραφε τῆς τε τρίποδος πέρι καὶ πεντέποδος ἀποφαίνων ὅτι μήκει οὐ σύμμετροι τῇ ποδιαίᾳ. καὶ οὕτω κατὰ μίαν ἑκάστην προαιρούμενος μέχρι τῆς ἑπτακαιδεκάποδος. ἐν δὲ ταύτῃ πως ἐνέσχετο. |

Some remarks on the translation.

1) The meaning of the verb 'ἔγραφε' in 147d2 had generated a lot of discussions. The first meaning of '*graphein*', in ancient Greek, is 'to draw' and especially in mathematics 'to draw a figure'. The term 'ἔγραφε' must be translated "he was drawing". Many examples can be cited in Plato that tend in this direction (for instance, *Republic* II 365c, 373a, 378c; VI 493c ; VII 523c, 529b, IX 583c, 589b (bis), X 602b; *Timaeus* 19b6, 71c4; *Laws* II 656e, 669a, VI 769b (bis), 769c, 889b). Theodorus is a geometer, and we must consider that he is drawing figures, probably in the sand. Later, Theodorus will proudly claim he left 'words' for 'geometry' (165a) which here could only mean the 'demonstrations by drawings'. His close association with Protagoras, and the spontaneous admission by his pupil Theaetetus that 'science' [37] (thus including mathematics) is 'perception' provide a supplementary hint about the kind of geometry he is doing.

   Now, it has a second meaning: 'to write'. But nothing is said about any kind of support for it. This is also the case in Socrates' lesson in the *Meno*. Moreover, what would it mean to 'write' some results? The most natural situation, which is the one considered by ancient and modern scholars, is the use of a stick on the dust of the ground or in sand. Neither would require the use of any specific instrument, of which the text says nothing. While it is not difficult to trace some drawings on the sand or even, for very rough ones, on the dust, to try to write a demonstration on it would be totally unreasonable. [38] So we can safely consider that Theodorus' lesson was centered on some graphics drawn on the sand of the stadium. In other words, the style of Theodorus' proofs is essentially graphical. [39] His tools are the usual ones that we can reasonably expect of an ancient Greek geometer: a rule – i.e. a stick – and a compass – i.e. probably a cord with two sticks at its ends (cf. *infra*, §B.ii)). And the only property used by Theodorus is (see below) a simple consequence of Pythagoras' theorem (cf. *infra*, **Appendix I**), which had been established much earlier in Greek mathematics.

2) There is a problem about the 'foot' used as unit (τῆς τε τρίποδος πέρι καὶ πεντέποδος). Actually, a foot is a unit of length, but Theaetetus employs it here as a unit of area. In the next sentence, it is said that Theodorus showed the incommensurability of these 'powers' of 3 and then 5 'feet'. Since 3 and 5 are commensurable to 1 (cf. *infra*, **Appendix III)**, the 'powers' refer necessarily to the sides of some areas of 3, 5 'feet'. Thus, the 'foot' is a unit of area and can be interpreted as a shortcut for 'a square of side 1 foot'. [40] Thus, when it is used as such, we will write it between quotes. As a matter of fact, many translators did not hesitate to translate '*podiaios*' here by 'foot-square'.

3) Against Burnet, we maintain ἀποφαίνων, which is the reading of all the manuscripts save one. [41]

4) Note the singular μήκει which refers to a plural; the same phenomenon occurs again below. This is a singular of totality.

5) The specification κατὰ μίαν ἑκάστην προαιρούμενος indicates that Theodorus is dealing with all these cases in a single lesson. The importance of this remark will appear later on (cf. §IV.C.ii).

6) On μέχρι τῆς ἑπτακαιδεκάποδος, there has been a polemic between Burnyeat and Knorr over whether 17 was included among the whole numbers studied by Theodorus. As we shall see, this controversy is as sterile as the one about *the* meaning of 'powers', and disappears in our interpretation.

### IV. Theodorus' lesson

Theaetetus' account of Theodorus' lesson (147d1-6), can be divided in:
   a) The construction of 'powers' (147d1-3)
   b) The proofs of the incommensurability of 'powers' (147d3-6).

### A. Theodorus' sequence

As odd as it may seem, there has been no real discussions on the integers considered by Theodorus, for the answer seemed obvious. However Theaetetus does not give an account of them. Instead, he uses an 'abbreviation': *3, 5, …, 17*. What we know is that this sequence begins with *3*, the next term is 5, and it stops at *17*, but we do not know which integers are considered between *5* and *17*. Moreover, there have been heated debates about the last integer *17*, whether or not it needs to be included among those studied by Theodorus, and whether he said anything about it. Historians of mathematics and commentators are deeply divided on the question. [42]

Conversely, since the Antiquity, commentators have wondered why Theodorus begins at *3* instead of *2*, for the square root of *2* is irrational [43] and its irrationality is easier to prove. [44] Most modern historians accept Hieronymus Zeuthen's explanation, [45] though sometimes with some hesitation: [46] Plato leaves it out because the irrationality of the square root of 2 was well-known.

They are certainly right that it was well-known in Theodorus' time, but this is not the reason for its absence from his lesson. Indeed, an important point for this question is almost universally neglected: the demonstration of the irrationality of the square root of *2* reduces the general problem of the irrationality of the square roots of integers to the odd integers alone. [47] As a direct consequence, Theodorus' sequence is easy to deduce: *3, 5, …, 17* is simply the natural sequence of the first odd integers (excluding the unit [48]) and ending at *17*.

According to the standard interpretations, Theodorus avoided perfect square integers. However, according to Theaetetus' account, this is not the case. True, translated into modern terms for brevity, the boy says that the mathematician showed that the square roots of *3* and *5* are irrational. But then he took 'each case in turn till he came to the one of seventeen' where he stopped. [49] Clearly, then, none was excluded. For such a natural sequence of integers, Theaetetus' shortcut was self-evident, and no further explanation was needed, either for Socrates or for Plato's reader.

There is also a historically important textual testimony that this sequence is the one actually studied by Theodorus. By chance, an anonymous commentary on the beginning of *Theaetetus* on a papyrus dated from the second century CE [50] covers the mathematical part entirely. The Commentator explains, among other possibilities, the absence of 2 in Theodorus' sequence, by the fact that the mathematician did not consider even integers, and this is the author's preferred explanation. At the very least, this shows the existence of a long tradition that claims that only the odd numbers were discussed in Theodorus' lesson. The reason was forgotten, for the entire theory was lost and replaced by the new theory of irrationality as found in Euclid's *Elements*.

B.  **The construction of 'powers'** (147d1-3)
i)  **Theodorus' drawings**
   a)  The preliminary drawing

In the very first sentence of the text (147d1), Theodorus proceeds by drawing some geometrical figures.

First, Theodorus traces two perpendicular straight lines of arbitrary lengths. In the course of the demonstration he will mark out on the horizontal line, the lengths of 1 foot, 2 feet, etc., as far as 9 feet. [51] The reason he does not need to go farther than 9 feet will be clear from the construction below. On the vertical axis, he will construct the lengths √3 feet, √5 feet √7 feet etc., as far as √17 feet. [52] He proceeds as follows.

   b)  The first drawing

He traces a horizontal line PZ of length of around 4 feet, marking each foot on it. Let P, Q, A, B, C be called these five points. Then let QT be the orthogonal straight half-line passing through Q, of length around 5 feet:

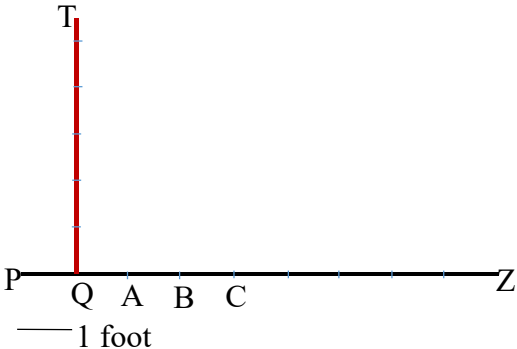

**Figure 1**

On this first figure, Theodorus draws a semi-circle of center A and or radius 2 feet, so that it cuts the line PZ at P and C. Let U be its intersection with the half-line QT (Figure 2):

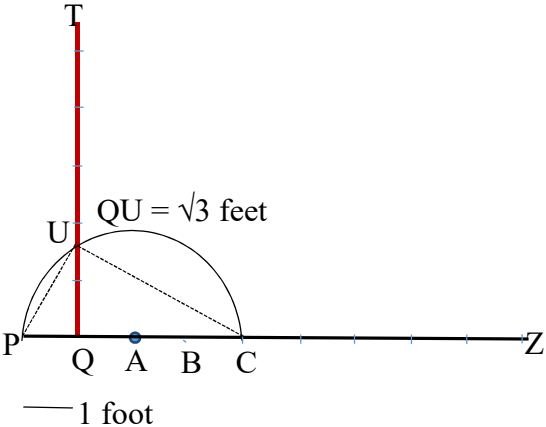

**Figure 2**

Then he draws the straight lines PU and CU, so that, since the angle CUP is inscribed within a semi-circle, it is a right angle. [53] Hence the triangle CUP is right-angled, so that from an easy corollary of Pythagoras' theorem (see **Appendix I**), we have:

the area of the square of side QU is equal to the area of the rectangle of sides PQ and QC. Since PQ = 1 foot and QC = 3 feet, we obtain: the area of the square of side QU is equal to 3 'feet', i.e. in modern writing QU = √3 feet. [54]

Let us also emphasize that there is no need to draw the lines PU and UC. This was done here only to show that the angle CUP is right, in order to apply the result of **Appendix I**.

    c)  The second drawing

On the previous drawing, he extends the line PC by two feet, and marks the successive points D and E (i.e. CD = DE = 1 foot). Then he draws the semi-circle of center B and radius 3 feet (so that it contains the points P and E). Let V be its intersection with the straight line QT (Figure 3):

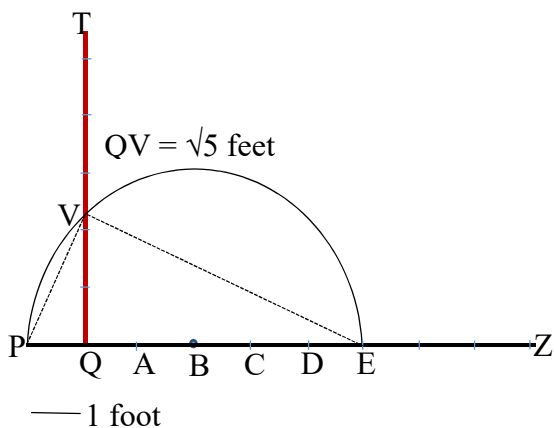

**Figure 3**

Once again, since the angle EVP intersects a semi-circle, it is a right angle, so that the corollary of Pythagoras' theorem gives:

The square on QV is equal (as area) to the rectangle on PQ and QE. [55] Since PQ = 1 foot and QE = 5 feet, we get: the area of the square of side QV is equal to 5 'feet', thus QV = √5 feet. [56]

Then he continues for 7, 9, 11, 13, 15 and 17 feet, 'in turn' ('κατὰ μίαν ἑκάστην'). Let us see what happens in the last construction.

d) The last drawing

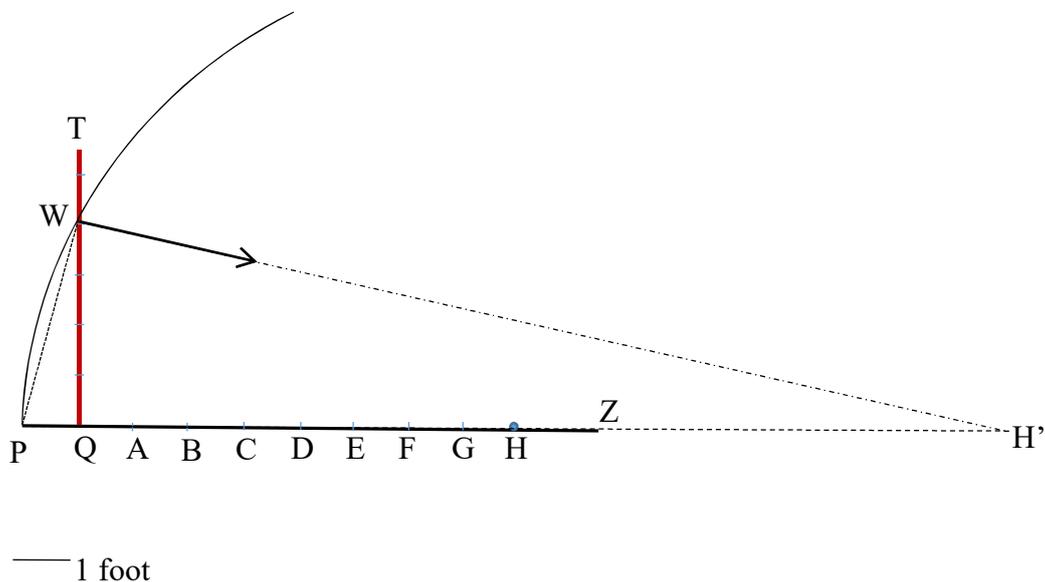

**Figure 4**

Let H be a point on the line PZ such that QH = 8 feet, and W the intersection of QT with the circle of center H and radius HP = 9 feet. Let us call H' (not necessarily drawn) the point on the extended straight line PZ such that HH' = 9 feet. As previously, H'WP is a right angle, and we have:

The square on QW is equal (in area) to the rectangle on PQ and QH'. [57] Since PQ = 1 foot and QH' = QH + HH' = 8 feet + 9 feet = 17 feet, we get: the area of the square of side QW is 17 'feet', hence QW = √17 feet. [58]

There is an historical testimony to Theodorus' construction as given above. It is found once again in the *Anonymous commentary*. [59] The drawing on the papyrus, illustrating the construction in a particular metrical case, is essentially that of proposition 14 of Book II of the *Elements* (with its accompanying drawing). [60]

ii) **Scaling or no scaling?**

A question now arises: since Theodorus uses a specific unit of length, the foot, did he really make the drawing according to this length or did he consider a 'symbolic' foot? In other words, did he consider a foot as a unit different from the ancient Greek unit of length, so that he might draw some smaller figure? The answer to this question does matter. First, some authors have even appealed to the reason that Theodorus could not go beyond 17 feet because of the dimensions of the stadium. [61] Others have studied whether early Greek geometers used the foot as a measure or simply as a symbolic unit. [62] As argued at the beginning of this article, one fundamental constraint in our interpretation of Plato's text is to take his words at face value. If in such a petty point the description is realistic, we could reasonably infer that Plato took care of the smallest details when he wrote a scene, because, contrary to the poets, [63] he knew about what he described either directly or from some specialists. For instance, he certainly attended lessons on the subject of irrationality, either by Theodorus or some other geometers.

Our purpose is to show, not only that a foot was indeed 1-foot long, [64] but that this is required in any real lesson, for the opposite would lead to absurdities.

First, there is no sense in speaking about a foot if any unit gives the same result. Moreover, the absence of any mention of scaling in the text confirms that a foot was really to be drawn as a foot. So, let us inquire whether this is possible, or whether it is just a product of Plato's imagination, or even a mistake of his part. We emphasized previously that Theodorus did not need to draw either the line WH' or the complete circle PWH' (i.e. of center H and radius 9 feet). [65] What he needed was the point W (cf. *supra*, Figure 4), which gave him the last side of the square he drew for his audience. In order to obtain it, he needed only to draw an arc of this circle around the line TQ, hence a straight line of length 9 feet, i.e. approximately 2.7 meters (to get P, Q, A, B, … and H) and small arcs of circles of radii at most 9 feet.

The final drawing on the ground (in the sand) should have more or less the following form:

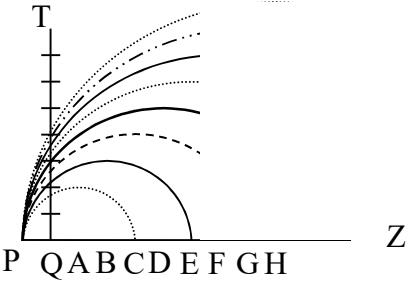

**Figure 5**

Thus, Theodorus drew eight successive arcs of circles, corresponding to all the odd integers between 3 and 17. What kind of tools did he need? a stick, a cord at least 9 feet long, with a knot at every foot, and two small sticks at its ends to draw the different arcs of circles of radius 2, 3, …9 feet, and perhaps a *gnomon* [66] to obtain better right-angled lines. Since all of this is so easily done with such simple equipment, there is no need to hypothesize any reduction in the size of the figures.

*Conversely*, there is one main reason **not** to proceed through scaling. We may certainly assume that his pupils were standing around him while he was tracing the successive drawings on the ground. Hence, Theodorus needed to get figures large enough so that the boys in his audience see them. Moreover, it was essential to be able to differentiate between the successive lengths – i.e., in modern language, √3, √5, …, √17 feet – the different objects of a lesson made 'each case in turn'. Since the last members of the sequence are √15 and √17 feet, their difference is around 0.25 feet. Hence, this is also the distance between the intersections of the line QT with the two last circles. With a scale of 1:1, it is around 7.5 cm, certainly close to the minimum needed to distinguish the intersections of two arcs of a circle roughly drawn in the sand (cf. Figure 5 above). Let us also remark that all these simple drawings are easy enough to be carried out in the time of a reasonably long lesson, leaving plenty of time for the proof of incommensurability.

We can safely conclude that Theodorus, or any geometer giving a similar lesson, would have drawn the sequence of lines of √3, √5, …, √17 feet (the intersection of the line QT with the different circles of centers A, B, …, and H) such that one foot would be close to 30 cm long.

### C. The proof of the incommensurability of 'powers' (147d3-6)
### i) The position of the problem

a) The geometer then proceeds to the proofs (147d2-4). The first result of incommensurability was the incommensurability, in Greek mathematics, of the diagonal to the side of a square. But this proof provided something more. The whole problem was reduced to the odd integers alone i.e. to finding a solution to the following problem: for any odd [67] integer $n$, whether there exist two odd integers [68] $p$ and $q$ such that: $\sqrt{n} = p/q$ which entails that $n = p^2/q^2$ or $nq^2 = p^2$. [69] In terms closer to Greek mathematics, this means that the side of the square of area $n$ is equal to the ratio of $p$ to $q$, which entails that $n$ is equal to the ratio of the squares of sides $p$ and $q$, which in turn entails that the square of side $p$ is equal to $n$ times the square of side $q$. In the first case (the existence of $p, q$) $\sqrt{n}$ is commensurable to the unit, otherwise it is not. Therefore, Theodorus is concerned only with odd integers.

b) Plato does not state Theodorus' method, not because it is useless for the understanding of Theaetetus' account, as claimed by M. Burnyeat, [70] but because it was well known to his readers. Besides, as far as we know, the proof given below is the only one that is both consistent with Theaetetus' account and based on results known by early Greek mathematicians. The main tool is the elementary 'remainder result', stating that the remainder after division by 8 of an odd square integer equals 1 (or in more scholarly terms, it is equal to *1 modulo 8*). [71] This proof is moreover closely connected with Theodorus' sequence defined in the previous paragraph.

c) Historians of mathematics who wished to account for Theodorus' result made no reference to the 'remainder result' prior to Jean Itard, who was, as far as we know, the first to show its importance for Theodorus' demonstration. [72] One of the scholars most sensitive to this question was W. Knorr in his famous work Knorr (1975). However, he did not take into account the importance of the proof of the incommensurability of the diagonal (in modern terms, the irrationality of the square root of 2) for this problem. Because he took the Platonic text seriously, he had noted that Plato imposed certain conditions upon any reconstruction of Theodorus' demonstration, and, as he rightly claimed, only his reconstruction verified them. His was thus the best attempt of his time (*ib.*, p. 96-97). However, it did not consider all the conditions necessitated by the text. In particular, the delicate condition of its duration: it had to be taught in a single lesson. (*ib.*, p. 193). This also implied certain translation difficulties, whence a polemic with M. Burnyeat concerning the verb '*enechesthai*'. [73] All this will be rectified in what follows.

ii) **On the supposed proof by '*anthyphairesis*'**

The most successful method among modern historians and mathematicians for explaining Theodorus' demonstration of incommensurability is the so-called proof by '*anthyphairesis*'. The term itself means 'alternate subtraction': given two integers, let them be $m$ and $n$, it consists in subtracting the lesser from the greater as many times as possible i.e.

if $m > n$, one writes: $m = kn + r$.

In modern terms, this means dividing the integer $m$ by the integer $n$; then $k$ is the quotient and $r$ the remainder of this division. E.g. for $m = 13$ and $n = 3$, we would get: $13 = 4 \times 3 + 1$ (since $5 \times 3 = 15$ is greater than 13), thus $k = 4$ and $r = 1$.

Since $n > r$ (otherwise it would be possible to subtract $n$ from $m$ once more) we could do the same for $n$ and $r$, and so on. [74]

The same process is possible when $m$ and $n$ are two arbitrary magnitudes, and a proposition of Euclid states that the process stops if and only if $m$ and $n$ are commensurable to each other. [75]

Many modern historians of mathematics have claimed that it was the method used by Theodorus to prove that the square roots of 3, 5, … up to 17 are incommensurable with the unit, by showing that the process is unlimited (i.e. that we get an infinite repetition of the same process). This process, of algebraic origin, is long and complicated, especially when translated into a geometrical approach. Moreover, it supposes a mastery of infinite processes not found in ancient Greek geometry. [76] However, the main problem is that it is inconsistent with Plato's text.

As a matter of fact, according to the text, the method needs to respect eight conditions (which include of course the ones required by Knorr):
   i) It must begin with 3, not with 2.
   ii) It is carried out 'each case in turn', not by means of a general result.
   iii) The fact that the process stops at 17 must be explained.
   iv) It must be consistent with mathematical knowledge at the time of Theaetetus' account.
   v) It must be able to be taught during one lesson for very young boys.
   vi) It must give rise to an erroneous generalization to all the integers.
   vii) It must be generalized to the case of cube roots of integers.
   viii) At the end of the proof, it must turn out that the quantity of 'powers' is unlimited.

The last three conditions are coming from the second part that will be studied in Brisson-Ofman (to appear).

The method using '*anthyphairesis*' verifies none of these conditions:

- Since the case of 2 is the simplest one, [77] Theodorus should have begun with this case in his lesson, and not with 3. Indeed it would be the best example for the young boys to understand such a difficult method, especially since the irrationality of the square root of 2 had long been known.

- To prove 'each case in turn' could perhaps be carried out during a year-long course, but certainly not in one lesson and that it would have made the boys flee from Theodorus, and possibly from mathematics.

- It certainly does not explain the stop at 17, since the argument that 19 is too complicated is not valid. On the one hand, as Knorr remarked, 19 is not much more difficult than the one for 13, and on the other if you are ready to carry out all the constructions needed for the proofs for all the integers between 3 and 17, you would certainly not be stopped by such a difficulty.

- There is no hint that such a method was available at the time of Theodorus, and indeed no effective use of an infinite sequence is found in ancient Greek geometry. On the other hand, the use of a finite algorithm of *anthyphairesis* to compute commensurable ratios is highly probable in early Greek geometry, but it is useless for the investigation of irrationality.

And from the second part of the passage, we get:

- There is no connection between such a proof and *Socrates*-Theaetetus' generalization as recounted at the end of Theaetetus' account.

- The property to show the irrationality of the square roots of integers through *anthyphairesis* (the infinite repetition of the same process), cannot be used for cube roots. Thus, the *anthyphairesis* cannot be generalized to the case of cube roots.

- There is no reason why the commensurability or incommensurability to the unit-foot of the side of the squares of area 3, 5, …, 15 (or even 17) times the unit-square through *anthyphairesis* would entail anything about the commensurability or incommensurability of other squares, and the proof has to be repeated for each different square. [78]

This shows beyond any reasonable doubt the inconsistency of this method with Plato's text. As noted, it also contradicts what is known of Greek mathematics in Antiquity. [79]

The success, in spite of this unlikeliness, of such a problematic method among mathematicians and historians of mathematics results from the flaws in the alternative methods proposed previously. They are essentially based on the spurious proposition X.117 in Euclid's *Elements,* which were already supposed to have given the first proof of incommensurability (cf. *infra*, **Appendix II**). In both cases, the demonstration gives a general answer to the problem of irrationality for all the square roots of integers studied by Theodorus, in flagrant contradiction with the text. In other words, we could say, like Socrates at the end of the *Theaetetus*, that it is, if possible, even more impossible ('ἀδυνατώτερον ἔτι ἐκείνων, εἰ οἷόν τε', 192b6) than by '*anthyphairesis*'.

iii) **The 'remainder result'** [80]

The method we present herein (*supra*, note 72) verifies all the above eight conditions. It is also interesting to compare its difficulty (for ancient geometers) to that of the method by '*anthyphairesis*'. [81] The present proof uses the simplest results and probably the oldest ones in the *Elements* (to be compared, for instance, with the entire chapter VI in Knorr (1975)).

The main part of the proof is the construction of geometrical figures, the result appearing directly in the drawings below:

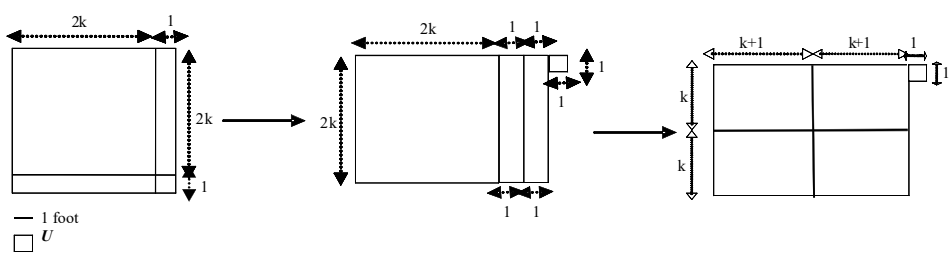

**Explanation**: We note by $U$ be the square-unit (i.e. of side 1 foot long). The square on the left is a graphical representation (using one foot as unit) of the square odd integer $n = (2k+1)^2$.

Moving the last raw (a rectangle of sides $(2k+1,1)$) to the right side of the square, we obtain a rectangle of sides $(2k+2, 2k)$ plus a small square of side 1. This is the middle figure.

The last figure is the same, but divided differently: four rectangles of sides $(k+1, k)$ plus the small square of side 1.

Now since either $k$ or $k+1$ is even, the last figure is equal to the sum of

a) 4 rectangles whose one side is even (hence a multiple of 8)

b) and a unit-square $U$.

Hence, the (area of the) square of side $(2k+1)$ is equal to (the area of) $U$ and a rectangle whose one side is 8.

From an arithmetic point of view, this means that the integer $n$ is the sum of an integer multiple of 8 plus 1. This entails [82] that the remainder of $n$ divided by 8 is equal to 1. [83]

**Remark 1**. However, the discovery of the result was probably done through the same method described by Theon of Smyrna (2[nd] century CE) to get **all** the square integers successively: adding the successive odd integers, beginning with 1 (*On Mathematics Useful for the Understanding of Plato*, XV). This is detailed in a figure, a little further on (*ib.*, XXV). In particular, all square odd integers are obtained in the following way:

1; 1 + (3+5) = 1 + **8** = 9 = $3^2$; 9 + (7+9) = 9 + **16** = 25 = $5^2$; 25 + (11 + 13) = 25 + **24** = 49 = $7^2$; 49 + (15 + 17) = 49 + **32** = 81 = $9^2$; etc.

Thus the difference of two successive odd square integers is always a multiple of 8, so that divided by 8, their remains are the same i.e. equal to 1. Q.E.D.

**Remark 2.** It is not difficult to write the above demonstration in the style of the arithmetical books of Euclid's *Elements*, and, in Remark 4, we will give an example of such rewriting for the next and longer proof.

It remains to be seen how this result concerning the integers alone, and not (what we call) 'rational numbers', may be used. This is not so obvious since all the results here concern integers, while Theodorus' purpose is to consider not whether the sides of the square of area 3 'feet', 5 'feet', …, 17 'feet' are integers or not, but whether they are or are not commensurable to the unit. [84] This explains why it was ignored so long.

### iv) Theodorus' characterization of quadratic commensurability

Let $n$ be an odd integer such tat the side of the $n$-foot square being commensurable to the 1-foot unit. According to §i).a), *supra*, there exist two odd integers $p$ and $q$ such that: $nq^2 = p^2$. Thus:

1) Since $q$ is odd, from the 'remainder result', we get that for some integer $k$:
$nq^2 = n (8k+1) = 8nk + n$.
Since $nq^2 = p^2$, and $p$ is odd, we may apply the same result to both $nq^2$ and $8nk + n$, so that: the remainder of $8nk + n$ divided by $8$ is equal to $1$.
Hence, the same holds true for $n$, [85] thus:
**the remainder of $n$ divided by 8 is equal to $1$**.                                                    (*).

2) Since obviously none of the integers 3, 5, 7, 11, 13, 15 verify this condition, the side of the squares of (area) 3, 5, 7, 11, 13, 15 'feet' are not commensurable to the 1-foot unit. And for the cases of 9 and 17 'feet', the 'remainder result' fails (cf. next paragraph).

**Remark 3**. The above algebraic proof is obviously anachronistic. Below we give a geometrical proof in which the result, as in the proof of the 'remainder result', is **read** from the drawings, justifying Theodorus' claim he left 'words' for geometry (i.e. for drawings). In **Appendix IV**, the interested reader will find the graphics that Theodorus may have drawn 'each in turn' for his prof.

### v) The demonstration

1) He began with the case of the 3-'foot'-square, i.e. with the following drawing:

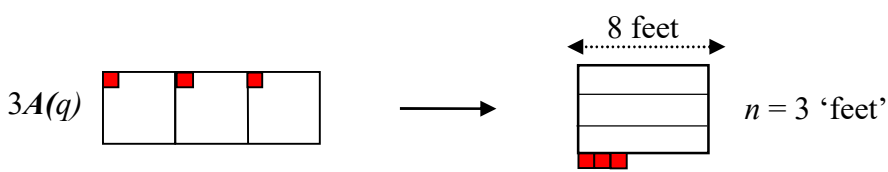

From the above drawing, for any odd integer *q*, we see that 3 times **any** square *A(q)* of side *q* is equal (in area) to a rectangle whose one side is 8-foot long, plus 3 times the square unit *U*.

But according to (\*), if the side *m* of the 3-foot square was commensurable to the unit of 1-foot long, there would be **an** odd integer *r* such that 3*A(r)* would be equal to a square *A(p)* for some odd integer *p*. Thus from §iii), *supra*, *A(p)*, thus 3*A(r)* would be equal to (the surface of) some rectangle whose one side is 8-foot long, plus **one** *U*. But the above figure shows that this is obviously false. Hence the side *m* is not commensurable to the 1-foot long side. Algebraically $\sqrt{3}$ is not rational.

2) The case of the 5-'foot' square. One has the following drawing:

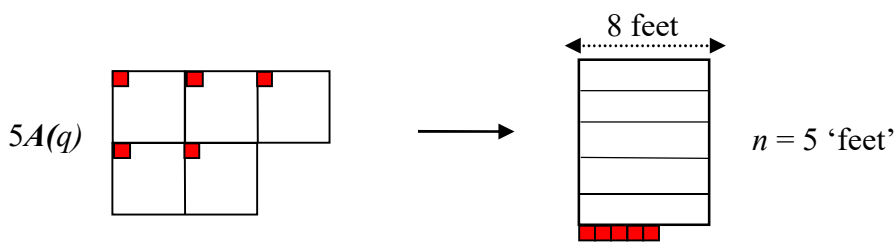

As previously, the drawing shows that for any odd integer *q*, 5 times any square of side *q* is equal (in area) to a rectangle whose one side is 8 feet long, plus 5 times the square unit *U*. But according to (\*), if the side *m* of the 5-foot square was commensurable to the 1-foot long side, there would be an odd integer *r* such that 5*A(r)* would be equal to a square *A(p)* for some odd integer *p*. The above figure shows this is impossible, hence the side *m* is not commensurable to the 1-foot long side. Algebraically $\sqrt{5}$ is not rational.

3) Then he would continue with 7-'foot' square. Its area is equal to a rectangle whose one side is equal to 8 feet long, plus 7 unit-squares *U*. Thus once again, its side *m* cannot be commensurable to the 1-foot long side, for as in the previous cases, it would need to be equal to some rectangle whose one side is equal to 8 feet plus one unit-square *U* (cf. Figure 8 in **Appendix IV**, *infra*). Algebraically $\sqrt{7}$ is not rational.

4) For the 9-'foot' square, it is equal (in area) to a rectangle whose one side is equal to 8 feet plus *U*. Thus, it verifies the necessary condition for its side to be commensurable to 1-foot (cf. Figure 8 of **Appendix IV**, *infra*). Actually, its side is 3-feet long, which is an integer, hence *a fortiori* commensurable to the 1-foot long side. Algebraically $\sqrt{9} = 3$ is rational.

5) Thus Theodorus 'goes on in this way, taking each case in turn', with the 11-'foot', 13-'foot' and 15-'foot' squares. For any odd integer *q*, 11 times (respectively 13 times, 15 times) *A(q)* is equal (in area) to some rectangle whose one side is 8 feet plus 3*U* (respectively 5*U* and 7*U*, cf. Figures 8 and 9 in **Appendix IV**, *infra*). Hence, as in the three first cases, the side of the squares of 11-'foot', 13-'foot' and 15-'foot' squares cannot be commensurable to the 1-foot long side. Algebraically: $\sqrt{11}$, $\sqrt{13}$ and $\sqrt{15}$ are not rational.

6) 'Till seventeen feet'. As shown in the last drawing of Figure 9 (**Appendix IV**, *infra*), for any odd integer *q*, 17*A(q)* is equal (in area) to a rectangle whose one side is 8 feet long, plus one unit-square *U*. Thus, as in the case of 9 feet, it verifies the necessary conditions for its side being commensurable to the 1-foot unit. But, contrary to the former, its side is obviously **not an integer**, for *17* is greater than *16 = $4^2$* and less than *25 = $5^2$*, and there is no integer between *4* and *5*. However, it is impossible to say anything about its commensurability or incommensurability. Hence, Theodorus stopped here and did not conclude, either because he wanted his students to find the result by themselves (as claimed by Árpád Szabó (1978)[86]), or, more probably, because he wanted them to attend the following lesson. [87]

**Remark 4.** It is easy to give the same demonstration more in line with the style of Euclid's *Elements*. We would have something like the following (we put our commentaries inside brackets):
Given an odd square [integer] A and U a unit; let U be subtracted from AB and CD the remainder.
Let EF be A added twice to itself.

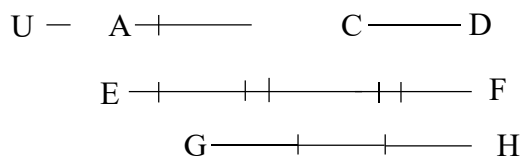

Figure 6

EF is also three times CD plus three units [evidence of the equality *3(x+1) = 3x + 3*]. Hence, since CD is measured by 8 [from the 'remainder result'], the same holds true for the remainder of three units subtracted from EF [since EF - 3 = 3CD]. Let it be GH. Thus, EF subtracted by one unit is not measured by 8, for otherwise the remainder of its difference with GH, that is 2, would be measured by 8 [evidence used implicitly for instance in the proof of proposition X.3] which is impossible Q.E.D.
Then it would go for 5A, 7A, …, 17A.

The above proof is much quicker and more elegant, but it is unlikely to have been Theodorus method:
- It is certainly less graphical.
- From a practical point of view, it would certainly not be easy to draw the last (very long) lines as far as 17A, especially since A cannot be too small a line.
- The use of a unit as the 'foot' is irrelevant.
- Last but not least, the boys would not have developed a graphical theory of integers, squares and rectangles, which would have been useless with this method (for this last argument, cf. Brisson-Ofman (to appear)).

**V. Conclusion**

This analysis of Theodorus' mathematical lesson shows, in a very concrete case, the extraordinary realism Plato displays when he accounts for events or depicts them. We have shown that, in the context of the mathematics of his own time, it was possible to give a lesson to young boys that respected the indications of the Platonic text down to the letter. It must be emphasized, of course – and this seems to have been too often forgotten – that he could hardly have done otherwise, since his characters were among those whom his readers knew, either by their speeches, or through the intermediary of parents or relations, or, finally, personally. In this context, the philosophers' friends and enemies both made any approach that was incompatible with the facts or the characters unacceptable. However, the choice of a work in dialogue form whose interlocutors were familiar to any Athenian citizen, was made by Plato who voluntarily does not hesitate to make use of ordinary, daily life, and hence leave himself open to criticism. In a forthcoming article, we study the second part of this passage, the work carried out by Theaetetus and his friend *Socrates*, on the basis of the lesson we have analyzed here. Unlike the usual interpretation, which considers it in isolation, or even ignores it, we will set forth its interactions with the rest of the dialogue, in which, moreover, mathematics are omnipresent, and we point out the importance of the "mathematical part" for the overall understanding of the dialogue as a whole.

**Appendix I. A corollary of Pythagoras' theorem**

Let $x$, $y$ and $z$ be three lengths. We note by $A(x)$ the area of the square of side $x$ and $R(y,z)$ the area of the rectangle of sides $y$ and $z$ (in modern symbolism: $A(x) = x^2$ and $R(y,z) = y \times z$).

**Result 1**

Given BDO be a right triangle, OB its hypotenuse and DH its height. Then

    i) $R$(OH,HB) = $A$(HD) i.e., HD is the mean proportional between OH and HB.
    ii) $R$(OH,OB) = $A$(OD) i.e., OD is the mean proportional between OH and OB.

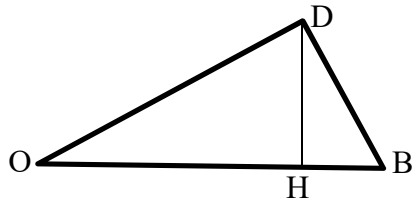

Figure 7

The translation into modern terms would be:
i) OH × HB = HD² or HB/HD = HD/OH and ii) OH × OB = OD² or OB/OD = OD/OH. [88]

**Appendix II. The incommensurability of the diagonal with respect to the side of a square**

1. The usual proof is based on the spurious proposition X.117 that is found in some manuscripts of Euclid's *Element* and widespread in the literature. Its simplicity (for moderns!) explains why it was so easily adopted, although it contradicts with the textual evidences. It is roughly the proof used nowadays. [89]
2. Because of its shortfalls, another method was proposed at the end of the 19th century by the Danish mathematician turned historian of mathematics Hieronymus Zeuthen, the 'anthyphairesis'. We have considered it in IV.ii), *supra*, and established its inconsistency with Plato's text. It was previously shown in Ofman (2010) that it could not be the original proof of incommensurability in Greek mathematics. In conclusion, both of the most popular methods for the origin of the incommensurability fail the test of the oldest extant texts.
3. In the same article we gave a method using old results, usually attributed to ancient Pythagorean or even earlier mathematics and, first and foremost belonging to the early theory of the 'odd and even' as claimed in Aristotle's *Analytics*.
4. The incommensurability of the diagonal to the side of the square means in modern terms that the square root of 2 is irrational (see **Appendix III**, *infra*). A simple consequence of the demonstration is that for any even integer $n$ there are two possibilities:

    - Either $n = 2^u \times v$ with $v$ and $u$ odd, in which case the square root of $n$ is irrational. A particular case is $n = 2$ ($u = v = 1$) but also $n = 6$ ($u = 1$, $v = 3$).

- Or $n = 2^u \times v$ with *v* odd and $u = 2s$, even. Thus, the square root of *n* is rational (respectively irrational) if and only if the square root of *v* is rational (respectively irrational). An example is $n = 12 = 2^2 \times 3$; its square root is equal to twice the square root of 3, hence it is rational/irrational if and only if the square root of 3 is rational/irrational ($\sqrt{12} = \sqrt{4 \times 3} = 2\sqrt{3}$). [90]

Therefore, the problem is **reduced** to the case of the odd integers. However, this does **not** mean that we know whether the square root of any even integer is or is not rational!

A natural question is why this method was not generalized, or at least used for instance for 3, since there are only 3 possibilities (instead of 2 for the previous case)? Even if we do not know whether or not it was tried in other cases, this is indeed unlikely. The main problem is with the so called 'Euclidean division', [91] which is called 'repeated subtraction' in Euclid's *Elements*. The question of the parity (or measure of parity in the sense of successive divisions by 2, for instance in Plato's *Parmenides*, 144a) is very special, to the point that all arithmetic is defined several times in Plato's texts as the 'knowledge of the odd and even' (for instance, *Gorgias*, 451b). Thus, any attempt to copy the previous process for integers greater than 2 needs to be considered not as a generalization (as for us) but as a completely different method.

**Appendix III. Some basic notions and correspondences between ancient Greek and modern mathematics**

We will first consider here some basic notions concerning mathematical incommensurability and irrationality.

For instance, while ancient Greek mathematicians did not speak of 'irrational **numbers**', they had some representation through geometrical 'irrational **magnitudes**'.

An '**irrational**' ('*alogos*') magnitude is defined on the basis of the above-mentioned incommensurable ratio. Thus let U be an arbitrary magnitude taken as a unit, e.g. a straight line. Any magnitude (e.g. a straight line) B is rational (a '*logos*') if the ratio of B to U is commensurable. If this ratio is incommensurable, B will be called irrational ('*alogos*'). Another term synonymous to '*logos*' ('**rational**'), the word '*rhêtos*'**,** was already used by Plato. [92]

Let us now compare the meanings of the terms '*logos*' ('**rational**') and '*symmetros*' ('commensurable').

 '**Commensurability/incommensurability**' vs '**rationality/irrationality**'**.** The former notions are relative ones and involve two objects (for instance two straight lines). The latter ones qualify a given object (for instance a straight line), but depend on an arbitrary object, the given unit.

- For instance, as seen in the previous paragraph, the diagonal and the sides of a square are **incommensurable** to each other, and this is independent of any other data.

- Now if we choose the side of a square as a unit, from the definition of **irrationality**, the length of its diagonal is **irrational**.
- But if we choose as a unit the diagonal itself, its length (i.e. 1) is of course **rational,** since its ratio to the unit is 1, an integer! [93]

To better understand the difference between commensurability/incommensurability and rationality/irrationality, let us consider an example written in modern language:

While $\sqrt{8}$ and $\sqrt{18}$ are **irrational** (cf. **Appendix II**.3), their ratio is:

$\sqrt{18}/\sqrt{8} = (3\sqrt{2})/(2\sqrt{2}) = 3/2$, hence $\sqrt{8}$ and $\sqrt{18}$ are **commensurable**.

In ancient Greek mathematics, one would say that while the sides of the squares of area 8 and 18 are **irrational**, the same sides are **commensurable** to each other.

Let us now summarize the different notions in ancient Greek mathematics and in modern mathematics:

| Ancient Greek mathematics | Modern mathematics |
|---|---|
| ἀριθμός: Integer<br>Field: arithmetic | Integral or whole (positive) number<br>Field: arithmetic |
| σύμμετρος λόγος**:**<br>Commensurable ratio<br>Field: arithmetic and geometry | Rational number<br><br>Field: arithmetic |
| ἀσύμμετρος λόγος**:**<br>Incommensurable ratio<br>Field: geometry | Real number<br><br>Field: algebra |
| σύμμετρον μέγεθος**:**<br>Commensurable magnitude<br>Field: arithmetic and geometry | Rational number<br><br>Field: arithmetic |
| ἀσύμμετρον μέγεθος**:**<br>Incommensurable magnitude<br>Field: geometry | Real number<br><br>Field: algebra |

Something new appears in the third and in the last lines of the above table, the 'real number'. In modern mathematics these 'numbers' are ubiquitous, although their definition is far from easy. Since there were no such numbers in Antiquity (they were defined at the end of the 19th century CE), ancient Greek mathematicians used lines. Straight lines replaced the modern 'real numbers', so that almost everything we do in modern algebra, analysis, and of course geometry, was then done through the representation of straight lines. Roughly speaking, where a modern mathematician may use 'numbers', ancient Greek mathematicians would use straight lines (among other 'magnitudes'). Indeed, throughout Antiquity and until modern times, 'geometer' was synonymous with 'mathematician'.

However, the situation is more complicated, because of the almost universal use of ratios in ancient Greek mathematics. For instance, irrationality and incommensurability were defined by means of such ratios. In modern mathematics a 'ratio' is simply the division of two 'numbers', thus another 'number', so that the 'ratios' simply disappeared from mathematics. The very meaning of 'ratio' was less clear in Greek mathematics, so that it was more difficult to define on them operations as simple as addition or multiplication.

To summarize, we may say that rational/irrational (λόγος/ἄλογος or ῥητός/ἄρρητος) qualifies magnitudes (such as straight lines), while commensurable/incommensurable (σύμμετρος/ἀσύμμετρος) qualifies ratios. Even when 'commensurable/incommensurable' magnitudes are used, [94] the meaning is that **their** ratios are commensurable/incommensurable to some other magnitudes. [95]

**Appendix IV. The geometrical application of the Remainder result**

In the diagram below, the colored portions represent one unit square *U* i.e. a square of a side one foot long, and the white parts represent areas that are multiples of *8*.

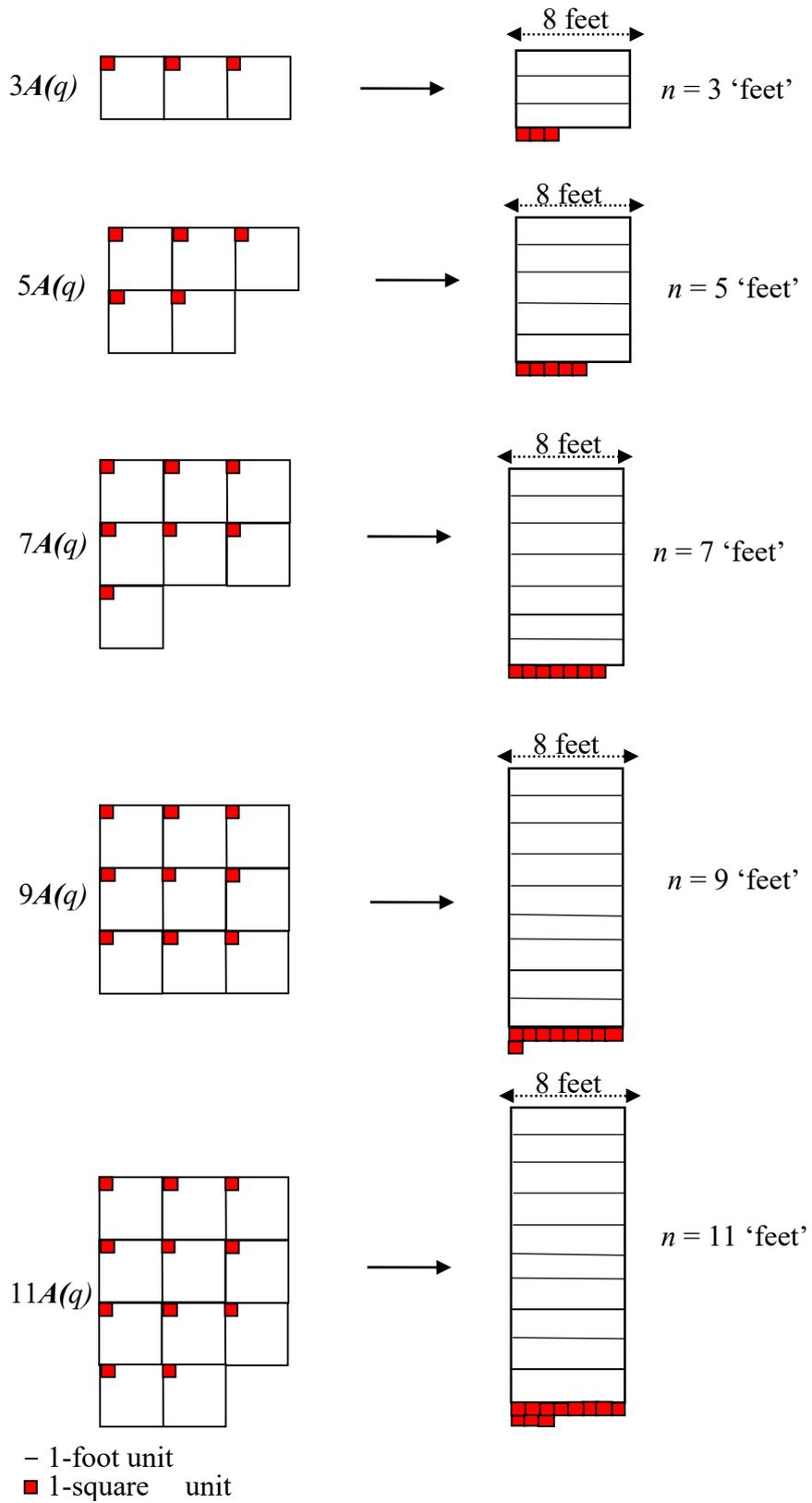

**Figure 8**

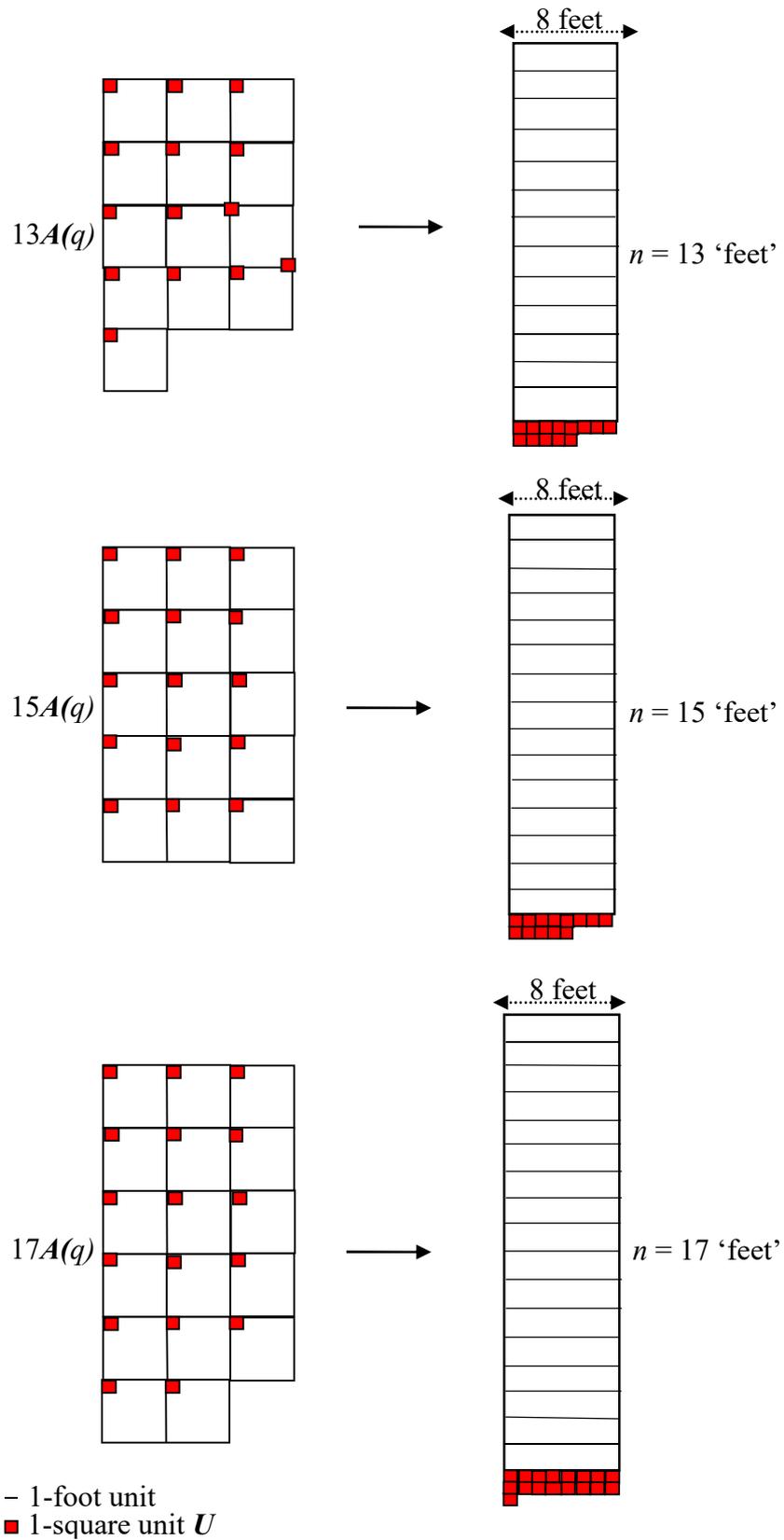

**Figure 9**

**Explanation.** On the left side, the squares represent an arbitrary odd square number $q^2$, and in the column its multiples by 3, 5, …17 are represented successively. The left to right arrow represents the passage from the left figure to the right one. It transforms truncated squares (i.e. squares minus one unit, hence 'gnomons') into rectangles, one side of which measures 8 feet long.

# NOTES

ignoring the controversy on the subject (for example, cf., Fowler (1999), p. 360. For a recent argued example against such a late dating, cf. Nails (2002), p. 275-278).

[23] See Brisson (2008).

[24] This does not mean, of course, that there were not important results in arithmetic (cf. for instance, book IX of Euclid's *Elements*). In particular, Knorr insists rightly that "Theodorus characteristically employed geometric constructions for the purposes of number theoretic studies" (Knorr 1975, n. 46, p. 101). As a matter of fact, such presentations by means of drawings lead to many complications (cf. reference in the note 25).

[25] 'So close was this association of theorem and diagram that the two terms might be used as synonyms.' (Knorr (1975), p. 72). See also A. Weidberg, *Plato's Philosophy of Mathematics*, 1955, p. 93; W. D. Ross, *Aristotle's Metaphysics*, 1924, I, p. 234, 295 and II, p. 268; J. L. Heiberg, *Mathematisches zu Aristoteles*, pp. 5-6. According to Eckhard Niebel, 'for Proclus a 'diagram' may stand for a 'theorem' as a **pars pro toto**' (Niebel (1954), p. 32 and 95, our emphasizes). Moreover, there are many recent works about the question of the proof and the use of the drawings by the ancient Greek geometers (see for instance Reviel Netz, *The Shaping of Deduction in Greek Mathematics*, Cambridge Univ. Press, 1999; Ken Saito, Traditions of the diagram, tradition of the text: A case study, *Synthese*, 2012, 186, p. 7–20).

[26] 'As for me [Theodorus], I bent away as early as possible from bare discourses towards geometry.' (ἡμεῖς δέ πως θᾶττον ἐκ τῶν ψιλῶν λόγων πρὸς τὴν γεωμετρίαν ἀπενεύσαμεν, 165a)

[27] Cf. also the long and detailed argument concerning this point made by A. Szabó in Szabó (1978), §3.1 (the 'proof' in Greek mathematics), p. 185-198.

[28] This is certainly not restricted to historians of mathematics. For example Burnyeat (1990): 'It is essential that we try to feel our way back into the Euclidean perspective, which is the closest we can get to the mathematics of Plato's time.' (p. 206). This leads to four pages on 'numbers' pitting against an alleged modern mathematical conception (i.e. Frege) over against what is alleged to be Plato's (i.e. Euclid's) conception.

[29] For instance, Netz (2004).

[30] He was confused for a long time throughout the Middle Ages with Euclid of Alexandria, the author of the *Elements*.

[31] *Supra*, note 22.

[32] At the end of the dialogue, Socrates says that he is going to see the examining magistrate, who will inform him of the charges brought against him by Meletus. We are therefore in 399 BCE.

[33] A Greek settlement located in the north-east of what is now Libya.

[34] Brisson (2000).

[35] That is through maieutics (Brisson (2008)). This is similar to the scene in the *Meno*, where Socrates made the young servant 'recall' the duplication of the square (82b-85b).

[36] Socrates immediately rejects this mixture, by changing Theaetetus' 'techniques' into 'sciences' (146d8).

[37] However, the remark about the meanings of 'epistêmê' at the end of §II.

[38] Heath's thesis that ''ἔγραφε' means 'to prove' is not supported by any testimony at the time of the dialogue (Heath (1921), I, p. 203). Moreover, even he admits the meaning of the sentence would be 'to prove through drawings'. For a detailed criticism of Heath's point of view, cf. Knorr (1975), p. 69-74.

[39] Though our analysis of Theodorus differs from Knorr's, we agree on this fundamental point: geometry, for him, is first and foremost graphical geometry. For a more detailed study on this question, Knorr (1975), section II of chap. III (p. 69-74, entitled 'on the role of diagrams'); also *supra*, note 25.

[40] The same form of 'shortcut' is used by Socrates in the lesson on doubling the cube in the *Meno* (also, *supra*, note 35).

[41] It 'significantly, appears in the dialogue-text reproduced in a papyrus commentary dating from the second or third century after Christ. Now, Burnet's reasons for bracketing 'ἀποφαίνων' are not clear. Furthermore, his edition was made several years before the papyrus commentary was published by Diels and Schubart. It is

certainly probable that had Burnet known of this commentary he would have retracted his bracketing of 'ἀποφαίνων'.' (Knorr (1975), p. 70).

[42] It is the subject of a heated dispute between W. Knorr and M. Burnyeat. The question was a philological one: what is the exact meaning of the verb 'ἐνέσχετο' in the sentence 'ἐν δὲ ταύτῃ πως ἐνέσχετο' in 147d6 (see the long note 88 in Burnyeat (1978), p. 512-513)? Does it mean that he necessarily had to stop **before** the last term (here 17 feet) because he met an obstacle there, as the former claims, or did he stop after it, at the next one (i.e. 19 feet)? It may seem a petty point, but certainly not according to the latter, who claims "It is not too much to say that large chunks of Knorr's impressive rewriting of the history of early Greek mathematics stand or fall by the thesis that this sentence means 'but in this one [sc. the 17-foot power] for some reason he encountered difficulty" (p. 513). Most scholars have followed Burnyeat's interpretation which is the standard one. As we will see later, according to our analysis, we agree with Knorr, against Burnyeat, that the proof given by Theodorus is essential for the understanding of Theaetetus' account, and that he did not consider the commensurability/incommensurability of the 17-foot power. Nevertheless, as we will see, Theodorus did draw the figure and studied this last case, but without concluding, leaving it open. Thus, the exact meaning of 'ἐνέσχετο' does not really matter.

[43] For the definitions of rational, irrational, commensurability, incommensurability, etc. and their relations, we refer to **Appendix III**.

[44] A proof consistent with textual testimonies, especially Aristotelian texts, is given in Ofman (2010). Although this proof is not absolutely necessary to understand Theodorus' method, it makes his lesson simpler and more natural (see also, *infra*, **Appendix II**).

[45] See for instance, Knorr (1975), p. 182.

[46] Taken as an explanation of last resort, because no others are available, for instance Burnyeat (1978), in particular p. 502-503.

[47] See Ofman (2010) and also **Appendix II**.4, *supra*.

[48] There is a discussion whether the unit was an integer ('arithmos') for the Greek mathematicians. The usual position, referring to a definition in Euclid's *Elements* (def. VII.2), is negative. However, many times in the *Elements*, it is treated as an integer. Moreover, in Plato's *Hippias*, the mathematician approves Socrates claiming 'one is odd' (*Great Hippias*, 302a). Last, in the second part of the passage, it needs to be considered as an integer according to Theaetetus' account (Brisson-Ofman (to appear)).

[49] 'καὶ οὕτω κατὰ μίαν ἑκάστην προαιρούμενος μέχρι τῆς ἑπτακαιδεκάποδος ἐν δὲ ταύτῃ πως ἐνέσχετο'. According to the classical interpretation, all the integers between 3 and 17 were studied except *4*, *9* and *16*. According to our interpretation, the only perfect square of the sequence is *9*, and is included in the sequence. This by the way solves what Jean Itard consider as a 'serious difficulty' ('une grave difficulté') concerning Theaetetus' account (Itard (1961), p. 34).

[50] A Berlin papyrus edited by Diels-Schubart (1905) under the title: *Anonymer Kommentatur zur Platonis Theaetet*, reedited with an Italian translation and commentaries by Bastiani-Sedley (1995).

[51] The foot (πούς) was a unit of length around 0.3 meters.

[52] Using modern symbolism here, the side of the square of area *n* is √*n* (for the question of the units, *supra*, III.1)).

[53] According to Diogenes Laertius (I 24), quoting Pamphilus of Epidaurus, this is one of the earliest results in ancient Greek mathematics. It is used by Aristotle as a trivial result for an example. It is finally the proposition 31 in book III of Euclid's *Elements*.

[54] Around 0.52 meter. For an explanation of 'foot' (inside quotes) considered as unit of area, see *supra*, III.1).

[55] In modern terms: $QV^2 = PQ \times QE$.

[56] Around 0.65 meter.

[57] In modern terms: $QW^2 = PQ \times QH'$.

[58] Around round 1.25 meter.

[59] Bastiani-Sedley (1995).

[60] It is written at the bottom of column XXXI of the papyrus (Bastiani-Sedley (1995), p. 346). It gives in the particular case of a rectangle of (area of) 3 'feet', the solution for constructing a square of the same area, using almost the same drawing found in Euclid's *Elements* (the one edited by J.L. Heiberg (1883-1885)):

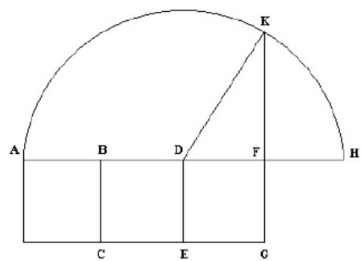

[61] Michael Paiow, Die Mathematische Theaetetsstelle, *Arch. Hist. Exact Sci.*, 27, 1, 1982, p. 87-99.

[62] Fabio Acerbi, In what proof would a geometer use the ΠΟΔΙΑΙΑ?, *Classical Quarterly*, 2008, p. 120-128.

[63] Cf. *Ion*, 536e-542a, in particular 541a-c.

[64] See supra, III.1).

[65] *Supra*, i).d).

[66] A sort of set square used in many fields, in particular for astronomical calculations. In Greek geometry, it is also the part of a square left when a smaller one has been removed from its corner.

[67] For the reason why it is not necessary to take even integers into account, see **Appendix II**.**4**, and for more details, Ofman (2010), §IV, p. 117-118.

[68] That the integers can be chosen to be odd does not arise from the definition; it is entailed by the proof of incommensurability (see *op. cit.*, previous note).

[69] We have the following inferences (as a matter of fact, equivalences here): $\sqrt{n} = p/q \Rightarrow n = p^2/q^2 \Rightarrow nq^2 = p^2$.

[70] Burnyeat (1978), p. 505.

[71] See for instance, Plutarch, *Platonic Questions*, II, 24, 1003f; also Iamblichus (*In Nicomachi*, ed. H. Pistelli, Leipzig, 1894, p. 90, 8), Diophantus (*Arithmetica* IV, 36). After a long analysis of the texts, Knorr's concludes this result needed to belong to 'the fundamentals of early Pythagorean arithmetic' (Knorr (1975), p. 154).

[72] Itard (1961), p. 34-36. Several historians followed Itard's suggestion, e.g. M. Caveing (Caveing (1998), p. 138-139), G.G. Granger (Granger (1998), p. 42), W. Knorr (Knorr (1975), p. 153-154), J. Vuillemin (Vuillemin (2001), p. 117-118). Indeed, the mere contemplation of a table of odd numbers leads to this property, which once understood, is easy to prove even within very early mathematics. Hence it is possible, as J. Itard claims, that it was already known by ancient Egyptian or Babylonian calculators (ib., p. 34). See also Remark 1, §IV.C.iii), *supra*.

[73] Burnyeat (1978), note 88, p. 512-513.

[74] For more details, see Ofman (2017a).

[75] For more details, see Zeuthen (1886); Becker (1933); van der Waerden (1963), *Science Awakening*, tr. A. Dresden, Science Editions, 1963 ( in particular 175-179); Caveing (1998), vol. 3, in particular chap. 3; Knorr (1975), p. 255-261; Fowler (1979), in particular p. 822-829.

[76] For instance Knorr (1975), p. 36. For the details of the geometrical constructions in the general case (i.e. including irrational lines) and their difficulties, p.118-126.

[77] For instance Caveing (1998), vol. 3, p. 123-127; Fowler (1979), p. 823.

[78] As a matter of fact, when we know the side $n$ of a square is incommensurable to the unit, for any integer $k$, it is obvious the same is true in the case of a square of side $k$ times $n$ (in modern terms, the irrationality of $\sqrt{m}$ entails the irrationality of $k\sqrt{m}$. But it is not what Theaetetus was meaning (cf. also the study of Knorr (1975), p. 84-85).

[79] For instance Knorr (1975)'s analysis in chapter 2, p. 29-35.

[80] A first version of this part was published in Ofman (2014).

[81] To get an idea of the complexity of obtaining Theodorus' results by 'anthyphairesis', cf. Knorr (1975), p. 118-126; cf. also for the case 7 feet, Vuillemin (2001), p. 147-150.

[82] We use the trivial property that for any integer, the remainders of $n$ and $n + 8$ after division by 8 are equal, thus the same holds true for $n + 8k$ for any integer $k$ (i.e., $n$ and $n+8k$ divided by 8 have the same remainder).

[83] In modern algebraic notation: $n = (2k+1)^2 = 4k^2 + 4k + 1 = 4k(k+1) + 1$. Now since either $k$ or $(k+1)$ is even, then $k(k+1)$ is even, so that $4k(k+1)$ is a multiple of 8. Hence, the remainder of $n = (2k+1)^2$ divided by $8$ is $1$ (cf. preceding note).

[84] In modern terms, the question is not whether $\sqrt{3}, \sqrt{5}, \ldots, \sqrt{17}$ are or are not integers, which is trivial, but whether they are or are not rational.

[85] Cf. *supra*, note 82.

[86] Szabó (1978), p. 79.

[87] Teachers were paid according to their audience: when Socrates is speaking with the so-called 'sophists', money matters are never far away. For instance, at the very beginning of the dialogue with Theodorus, the first subject is money, more precisely Theaetetus' money (140c6-d3). And further on in the dialogue, when discussing some theses of Protagoras, a friend of Theodorus, Socrates recalls that the former demanded 'huge fees' from his followers (161d).

[88] For an elementary graphical proof, Ofman (2017b).

[89] For its criticism, for instance, Ofman (2010), p. 87-90.

[90] For the details, *ib.*, p. 117-8.

[91] I.e. the division of an integer by another with dividend and remainder. For instance, 13 divided 'Euclideanly' by 4 is equal to 3 with remainder 1 ($13 = 4 \times 3 + 1$).

[92] Though Euclid's terminology in the *Elements* is somewhat different, it is the one largely used by the Greeks during Antiquity.

[93] This opposition, called by Pappus, 'by convention' vs. 'by nature', is already noted in his Commentary of Euclid's book X, §14 (Thompson (1930), p. 78 (p. 14-15)).

[94] The connection between commensurability and rationality was so well-known, that some authors would sometimes only explicitly mention one term of the relation, the other being supposed evident. For instance, Aristotle writes: 'the diagonal *is incommensurable*, since the contradictory proposition produces a false result.' (*Prior Analytics* I, 23, 41a29-30, transl. H. Tredennick, our emphasize). However, the meaning is clearly that the diagonal is incommensurable with respect to the side.

[95] For instance, the first definition of book V of Euclid's Elements: **Σύμμετρα μεγέθη** λέγεται τὰ τῷ αὐτῷ μέτρῳ μετρούμενα, … ('**Those magnitudes** are said to be **commensurable** which are measured by the same measure, …').